\documentclass[10pt]{amsart}

\usepackage{amsfonts,amsthm,latexsym,amsmath,amssymb,amscd, amsmath, mathrsfs,color}

\usepackage[cp1251]{inputenc}
\usepackage[russian,english]{babel}

\theoremstyle{plain}
\newtheorem{theorem}{Theorem}[section]

\theoremstyle{definition}
\newtheorem{definition}{Definition}[section]

\newtheorem{example}[theorem]{Example}

\begin{document}

\author[O.~Katkova]{Olga Katkova}

\address{Department of Science and Mathematics, Wheelock College, USA}
\email{olga.m.katkova@gmail.com }

\author[M.~Tyaglov]{Mikhail Tyaglov}

\address{School of Mathematical Sciences, Shanghai Jiao Tong University\\
and  Faculty of Mathematics, Far East Federal University}
\email{tyaglov@sjtu.edu.cn}

\author[A.~Vishnyakova]{Anna Vishnyakova}

\address{School of Mathematics and Computer Sciences, Kharkov National V.N.Karazin University}
\email{anna.m.vishnyakova@univer.kharkov.ua}

\title[Finite difference  operators]
{Linear finite difference operators with constant coefficients and distribution of zeros of polynomials}





\begin{abstract}

We study the effect of finite difference operators of finite order
on the distribution of zeros of polynomials and entire functions.

\end{abstract}

\maketitle


\setcounter{equation}{0}

\section{Introduction}\label{Section:intro}

One of the important problems in the theory of distribution of zeros of polynomials and transcendental entire functions
is to describe linear transformations which map polynomials having all zeros in a given region 
into the set of polynomial having all zeros in another given region. A very important case is that of both regions being 
equal to the real line.

\begin{definition}
\label{th:d1}
A real polynomial $P$ is called hyperbolic  (or real-rooted) if all 
zeros of $P$ are real or if  $P$ is identically zero.
\end{definition} 
As usual we denote by $\mathcal{HP}\subset \mathbb {R}[x]$ the set of  
hyperbolic  polynomials.

Hermite and, later,  Laguerre were, probably, the first  to study such type of problems systematically.
In 1914 P\'olya and Schur~\cite{polsch} completely described the operators acting diagonally on the standard monomial
basis $1$, $x$, $x^2$, \dots of $\mathbb{R}[x]$ and preserving the set of hyperbolic polynomials. Later the study of linear
transformations sending real-rooted polynomials to real-rooted polynomials was continued by many authors including
N.\,Obreschkov, S.\,Karlin, B.\,Levin, G.\,Csordas, T.\,Craven, K.\,de Boor, R.\,Varga, A.\,Iserles, S.\,N{\o}rsett, E.\,Saff etc.
Among recent authors it is especially worth to  mention P. Br\"{a}nd\'{e}n and J. Borcea~\cite{BrandenBorcea} (see also~\cite{BrandenBorcea2,BrandenBorcea3}), who completely
characterized all linear operators preserving real-rootedness of real polynomials (and some other root location preservers).

A natural extension of
polynomials with real roots is the so-called Laguerre-P\'olya class.

\begin{definition}
\label{th:d2} A real entire function $f$ is said to be in the {\it
Laguerre-P\'olya class}, written $f \in \mathcal{L-P}$, if
\begin{equation}\label{e2}
 f(z) = c z^n e^{- a z^2+b z}\prod_{k=1}^\infty
\left(1-\frac {z}{x_k} \right)e^{\tfrac {z}{x_k}},
\end{equation}
where $c, b, x_k \in  \mathbb{R}$,  $x_k\ne 0$,  $a \geqslant 0$,
$n$ is a non-negative integer and $\sum\limits_{k=1}^\infty x_k^{-2} <
\infty$. The product in the right-hand side of~\eqref{e2} can be
finite or empty (in the latter case the product equals 1).
\end{definition} 

This class is essential in the theory of entire functions due to the fact that these
and only these functions are the uniform limits, on compact subsets of
$ \mathbb{C}$, of polynomials with only real zeros. For various properties and
characterizations of the
Laguerre-P\'olya class see, e.g. \cite[p. 100]{pol}, \cite{polsch}, \cite[Chapter VII]{lev}, 
\cite[pp. 42--47]{HW},  \cite[Kapitel II]{O} or  \cite{CV}.

G.\,P\'olya  obtained probably the first results on $\mathcal{L-P}$-preservation properties of a 
linear finite difference operator.
In~\cite{pol1} he established that if
$f \in \mathcal{L-P},$ then  $f(x+ih) + f(x-ih)\in \mathcal{L-P}$
for every $h \in \mathbb{R}$. N.G. de Bruijn observed that this fact can be refined as follows \\

{\bf Theorem A} (\cite[Theorem 8]{Br}).  { \it  For arbitrary $ h, \alpha \in \mathbb{R}$ the 
linear operator 
\begin{equation}
\label{e33}
B_{h, \alpha}(f) (x) := e^{i\alpha}f(x+ih) + e^{-i\alpha} f(x-ih)
\end{equation} 
preserves the class $\mathcal{L-P}$.}

Our object of study is linear finite difference operators with constant coefficients.

Let $T : \mathbb{C}[x] \to \mathbb{C}[x]$ be a linear finite
difference operator of the form
\begin{equation}
\label{e3}T(P) (x) = \sum_{j=l}^m a_j P(x-j\lambda), 
\end{equation}
where $l, m \in \mathbb{Z}, l < m $, $a_j
\in \mathbb{C}, l\leq j \leq m, a_l \neq 0, a_m \neq 0,$  
$ \lambda \in \mathbb{C}\setminus \{0\}.$\\

Our present work was inspired by the paper~\cite{BKS} of P.\,Br\"{a}nd\'{e}n, I.\,Krasikov and B.\,Shapiro, 
where they studied linear finite difference operators with polynomial coefficients and with a 
real shift $\lambda.$ The authors made an attempt to transfer the existing
theory of real-rootedness preservers to the basis of Pochhammer symbols and to develop a 
finite difference analogue of the P\'olya-Schur theory. In particular, in \cite{BKS} it was 
proved that  a linear operator 
of the form (\ref{e3}) with a real shift $\lambda$ preserves the set of hyperbolic polynomials 
if and only if at most one of coefficients $a_j(x)$ is nonzero, and $a_j(x)$ is hyperbolic for such a $j$.\\

In the present paper we study linear finite difference operators of the form (\ref{e3}) with 
constant coefficients and with a complex shift $\lambda.$
We need to express such a  finite difference operator in terms of the shift operator.
\begin{definition}
\label{th:d3}
For every  $\lambda\in \mathbb{C}$ define the shift operator: $S_\lambda : \mathbb{C}[x] \to \mathbb{C}[x] $ by
$$S_\lambda(P)(x) :=P(x-\lambda). $$
\end{definition}
Obviously we have 
\begin{equation}
\label{e4}T(P) (x) = \sum_{j=l}^m a_j S_\lambda^j(P)(x).
\end{equation}

We will consider the generating rational function of the operator $T$: 
\begin{equation}
\label{g1} Q (t) = \sum_{j=l}^m a_j t^j.
\end{equation}

We give a description of linear operators of the form  (\ref{e3}) with an arbitrary complex 
shift which preserve the 
set of hyperbolic polynomials. 

\begin{theorem}\label{th:mth1}

Linear operator $T$ of the form (\ref{e3}) preserves the set of 
hyperbolic polynomials if and only if the following 
conditions are satisfied:

1.  $ \mbox{Re} \   \lambda =0; $ 

2.  $l = - m;$

3. All roots of the generating function (\ref{g1}) belong to the unit circle 
$\{ z : |z| = 1 \}$;

4. $ a_{-m}\cdot a_{m}\  \in (0; + \infty). $
\end{theorem}

{\bf Remark.} The assertions 2, 3, 4 of the above theorem mean that the generating function (\ref{g1}) is of the form 
$$ Q(t)= C\prod_{k=1}^{2m}\left(e^{-i\theta_k/2}\sqrt{t}+e^{i\theta_k/2}\frac{1}{\sqrt{t}} \right),$$
where the numbers $C$ and $ \theta_k \   ( k=1, 2, \ldots,  2m)$  are real. Thus, theorem \ref{th:mth1} states 
that every linear operator of the form (\ref{e3}) that preserves the set of  hyperbolic polynomials is a 
composition of  linear operators of the form (\ref{e33}), that is the operators of
 the form  $e^{i\alpha}f(x+ih) + e^{-i\alpha} f(x-ih), \quad h, \alpha \in \mathbb{R}.$ \\

It turns out that the operators described in Theorem \ref{th:mth1} are also strip preservers.
\begin{theorem}\label{th:mth2}

Let $b>0$ be a given number. Linear operator $T$ of the form (\ref{e3}) preserves the set of 
complex polynomials having all zeros in the strip $\{ z:\  |\mbox{Im}\  z|\leq b \}$ if and only if conditions 1-3 
are valid.

\end{theorem}

We note that the sufficiency of such  conditions for linear operator of the form (\ref{e3}) to preserve the set of 
polynomials with all zeros in the strip (and other interesting properties of such operators) was proved in \cite{Br}.

The fact that every linear operator of the form (\ref{e3}) that preserves the set of  hyperbolic polynomials is a 
composition of  linear operators of the form (\ref{e33}) motivates us to study such kind of operators  in more detail. Further it will be
more convenient for us to put $\alpha = \theta - \pi/2$ and to study the following equivalent form of the operator
\begin {equation} \label{e5} 
T_{\theta, h}(P)(x)=\frac{e^{i\theta}P(x+ih) - e^{-i\theta}P(x-ih)}{i}, \quad h>0,\   \theta \in\mathbb{R}.
\end{equation}

We note that in \cite{KTV} the complete description of all finite difference operators of the form $
\Delta(f)(z)=M_1(z) f(z+ih) + M_2(z) f(z-ih)$ (where $M_1$ and $M_2$ are 
some complex functions, $h>0$), preserving the Laguerre-P\'olya class was obtained.
\\

The following example is  important in the sequel. 
\begin{example} \label{E1} For $n\in\mathbb{N}$ we consider $L_n(x) =x^n \in \mathcal{HP}.$ It is easy to calculate that
\begin {equation} \label{E2} Q_n(x,\theta)= T_{\theta, 1}(x^n) =\frac{e^{i\theta}(x+i)^n-e^{-i\theta}(x-i)^n}{i}=
2 \sin\theta \prod_{k=1}^{n}\left (x-cot\frac{-\theta +\pi k}{n} \right ), \mbox{if}\  sin\ \theta \neq 0,\ \end{equation}
and for $\theta $ with $sin\ \theta =0$
\begin {equation} \label{E3} Q_n(x,2\pi m)=- Q_n(x,\pi+2\pi m)=Q_n(x,0)=\frac{(x+i)^n-(x-i)^n}{i}=
2 n\prod_{k=1}^{n-1}\left (x-cot\frac{\pi k}{n} \right ), \ m\in\mathbb{Z}. \end{equation}
\end{example}

We will denote by
\begin {equation}
x_k=x_k(\theta)=cot\frac{-\theta +\pi k}{n}, \ 
\ k=1,2,\cdots, N,\label{a3} \end{equation} the zeros of the polynomial $Q_n(x,\theta),$ where $N=n$ if $sin\ \theta \neq 0$ 
and $N=n-1$ if $sin\ \theta=0.$

We observe that all zeros of $Q_n(x,\theta)$ are real and simple. It is easy to show that for every hyperbolic 
polynomial $P$ all roots of  $T_{\theta, h}(x)$ are simple. We observe also that the minimal distance between
different zeros of $Q_n$ tends to zero when $n$ tends to infinity. So we can not use the limiting reasoning
to conclude that all roots of  $T_{\theta, h}(f)$ are simple for all $f \in \mathcal{L-P}.$ We proved the following theorem. 

\begin{theorem}\label{th:mth3}
For every $h>0,\  \theta \in \mathbb{R},$ and  every $f \in \mathcal{L-P},$  all the zeros of $T_{\theta, h}(f)$ are real and simple.
\end{theorem}

 For every hyperbolic polynomial $P$ we  obtained the estimation for the maximal and minimal roots of the image 
 $T_{\theta, h}(P).$ Let's denote by $\lambda(P)$ the {\it maximal} root of a hyperbolic polynomial $P$ and by 
 $\mu(P)$ its {\it minimal} root. We prove the following statement.
 
 \begin{theorem}\label{th:mth5}
\it For every  $P \in \mathcal{HP}$ , $\deg P =n \geq 1,$  $\theta\in\mathbb{R},  $ 
and every $h>0,$ we have 

$$\lambda( T_{\theta, h}(P)) \leq \lambda(P)+h\cdot \lambda(Q_n) \ \ \mbox{and} \ \ 
\mu( T_{\theta, h}(P)) \geq \mu(P)+h\cdot \mu(Q_n),$$
where polynomials $Q_{n}$ are taken from Example \ref{E1}.

 \end{theorem}

To formulate our next result we need further the following frequently used measure of zero separation 
for hyperbolic  polynomials.
\begin{definition}
\label{th:d4}
Given a  polynomial $P\in \mathcal{HP}, \  \deg P \ge
2,$ denote by $\mathrm{mesh} (P)$ the minimal distance between its
roots: $$ \mathrm{mesh} (P) := \min\limits_{1\leq j\leq n-1} (x_{j+1}-x_j)$$
for $ P=C(x-x_1)(x- x_2) \cdot \ldots \cdot (x- x_n),$  where  $x_1 \leq x_2 \leq \ldots \leq x_n.$ 
(If $P$ has a double real root, then $\mathrm{mesh}(P) =0$).  
\end{definition}

The following beautiful fact was discovered by M.Riesz in 1925 and probably initiated the 
study of mesh nondecreasing operators. \\

{\bf Theorem B} (M. Riesz, 1925). {\it Let $P\in {\mathcal HP},\  \deg P \geq 3.$ Then 
$\mathrm{mesh} (P^{\prime}) \geq \mathrm{mesh} (P) $. If all zeros of $P$ are simple, then 
$\mathrm{mesh} (P^{\prime}) > \mathrm{mesh} (P) $. }\\

An elementary proof of this theorem was given by A. Stoyanoff (\cite{Sto}).  It turns out that 
there are other linear operators that do not decrease the mesh of hyperbolic polynomials. 
The following result shows that if a hyperbolicity preserver commutes with the shift operators, it does 
not decrease mesh.\\

{\bf  Theorem C} (S. Fisk, \cite[p. 226, Lemma 8.25]{Fisk}).  {\it If  $A: {\mathcal HP} \to  {\mathcal HP}$ 
is a linear operator, and for all
$  b\in \mathbb{R}$ we have $A S_b = S_b A,$  then for every
 $ P\in {\mathcal HP}$ the following inequality holds:  $\mbox{mesh} (A(P)) \geq \mbox{mesh} (P).$} \\
 
 Note that S.Fisk formulated this theorem in other terms. It is not easy to recognize that 
 S.Fisk's theorem is the statement above. 

As we mentioned earlier in  \cite{BKS} it is proved that any nontrivial linear operator 
of the form (\ref{e3}) with a real shift $\lambda$ does not preserve the set of hyperbolic 
polynomials. But in  \cite{BKS} it is proved that a linear operator 
of the form (\ref{e3})  with a real shift $\lambda$ preserves the set of hyperbolic polynomials
having mesh not less than $\lambda$  if and only if all zeros of the generating rational function 
$Q(t) := \sum_{j=l}^m a_j t^{j}$ are real and non-negative.

 Since the linear operator $ T_{\theta, h}$ is a hyperbolicity preserver for every $\theta, h\in\mathbb{R}$ and 
$ T_{\theta, h}$  commutes with any shift operator,  it follows from Theorem C that $ T_{\theta, h}$ does not 
decrease mesh. We show that in the class of all hyperbolic polynomials of degree $n$ the polynomial $x^n$ is 
extremal in the following way.

\begin{theorem}\label{th:mth4}
\it For every  $P \in \mathcal{HP}$ , $\deg P =n \geq 2,$  every  $\theta\in\mathbb{R},\   \sin \theta \neq 0,$ 
and every $h>0,$ we have 
 $$\mbox {mesh}\  T_{\theta, h}(P)  \geq \mbox {mesh}\  T_{\theta, h}(x^n) ;$$  
For every $ \theta:  \sin \theta = 0,$ the statement of the theorem is also true for all $ n \geq 3.$
 \end{theorem}

In connection with theorems \ref{th:mth5} and \ref{th:mth4} the following natural question arises. 

{\bf Open problem.} {\it   To describe the image of the set of hyperbolic polynomials (of the set of hyperbolic 
polynomials of degrees not greater than a given $n$)  under the linear operator of the form (\ref{e5}).}\\
 
Our last theorem describes the asymptotic behavior of zeros of $T_{\theta, h}(P)$ for $h$ tends to infinity. 
Let $P_n(x)=x^n+ax^{n-1}+bx^{n-2}+\sum_{k=0}^{n-3}c_kx^k$ be a polynomial with complex coefficients.
For  $\theta \in \mathbb{R}, h>0$ consider the polynomial 
\begin {equation} D_n(x, \theta,h):=  T_{\theta, h}(P_n)(x)= \frac{e^{i\theta}P_n(x+ih)-e^{-i\theta}
P_n(x-ih)}{i}.\label{e6} 
\end{equation} 
The polynomial $D_n(x, \theta,h)$ has $n$ roots if $\sin \theta \ne 0,$ while it has only $n-1$ root if $\sin \theta = 0.$
Denote by $X_1(h,\theta), X_2(h,\theta),\ldots, X_{N}(h,\theta)$ the roots of this 
polynomial numerated
under the condition: $\mbox{Re}\  X_1(h,\theta) \leq \mbox{Re}\  X_2(h,\theta) \leq \ldots \leq \mbox{Re}
\  X_{N}(h,\theta),$  where $N=n$ if $sin\theta \neq 0$ and $N=n-1$ if $sin\theta=0.$ Our goal is to 
describe an asymptotic behavior of $X_j(h,\theta), \ \ j=1,2,\ldots,n,$ as $h\to\infty.$ 

\begin{theorem}\label{th:mth6} For every   $\theta\in\mathbb{R},$ 
and  $h>0,$   the j-th root of the polynomial $D_n(\theta,h)$ satisfies the asymptotic formula:
$$X_j(h,\theta)=x_j\cdot h-\frac{a}{n}+\left( \frac{a^2(n-1)}{2n^2}  - 
\frac{b}{n}\right)\frac{Q_{n-2}(x_j,\theta)}{Q_{n-1}(x_j,\theta)}\cdot \frac{1}{h} +
O\left(\frac{1}{h^2}\right), \quad h\to \infty ,$$
where polynomials $Q_{n-1}, Q_{n-2}$ and numbers $x_j$ are taken from Example \ref{E1}.
\end{theorem}

\setcounter{equation}{0}

\section{Proof of Theorems \ref{th:mth1} and \ref{th:mth2}  }\label{section:proof.thm1.1}

Suppose that a linear operator $T$ of the form (\ref{e3}) preserves the set of 
hyperbolic polynomials. For every  $n\in \mathbb{N}$ we consider a hyperbolic polynomial
$L_n (x) =x^n.$  We have
\begin{equation}
\label{f2} T(L_n) (x) = \sum_{j=l}^m a_j (x-j \lambda)^n = x^n 
\sum_{j=l}^m a_j \left(1-\frac{j \lambda}{x}\right)^n =:
x^n S_n(x) \in {\mathcal HP}.
\end{equation} 
Thus all the zeros of the rational function $S_n$ are real. 

Put $x=\frac{n}{y}, \  y\in
\mathbb{R}\setminus \{0\}.$ By our assumptions for every  
$n\in \mathbb{N}$ all the zeros of $S_n (\frac{n}{y})$ belong to 
$\{ z : \mbox{Im} \  z = 0 \}.$ The sequence $S_n (\frac{n}{y})$ converges 
uniformly on the compact sets to the entire function $$f(y) := \sum_{j=l}^m 
a_j e^{-j \lambda y}$$ as $\ n\to \infty.$ We conclude that all the
zeros of the entire function $f(y)= Q (e^{- \lambda y})$ are real. 

Let us find 
the zeros of $f$. We put $\lambda = \alpha + i \beta,\   \alpha, \beta \in \mathbb{R},$ 
and suppose that $z_0 \in \mathbb{C}\setminus \{0\}$ is a zero of $Q$. Then we 
solve the equation $$e^{-\lambda y} = z_0$$ and get $$y_k = - \frac{\log |z_0| +
i\arg z_0 + 2\pi k i }{\alpha + i \beta}, \ \  k\in \mathbb{Z}.$$
Thus $$\mbox{Im}\  (y_k) =  \frac{\beta \log |z_0| - \alpha \arg
z_0 - 2\pi k \alpha }{\alpha^2 +  \beta^2}, \  k\in
\mathbb{Z}.$$ Since by our assumptions $\mbox{Im}\  (y_k) = 0$ for
all $k\in
\mathbb{Z}$ we get $$\alpha =\mbox{Re}\  \lambda =0.$$ Whence $$\mbox{Im}\  (y_k) = - \frac{-
\log |z_0| }{  \beta} = 0, \  k\in \mathbb{Z},$$ and we obtain
that any non-zero root of $Q$ belongs to the circle  $\{ z : |z| = 
1\}.$ The necessity of the conditions 1 and 3 in theorem \ref{th:mth1} is proven.

 The necessity of the conditions 1 and 3 in the case of strip preservers can be shown analogously.  Suppose that a linear operator 
$T$ of the form (\ref{e3}) preserves the set of complex polynomials having all zeros in 
the strip $\Pi_b :=\{ z:\  |\mbox{Im}\  z|\leq b \}.$ Then from (\ref{f2})  we get that all the zeros of the rational 
function $$S_n(x) =\frac{T(L_n) (x)}{x^n}$$ belong to $\Pi_b.$

 Put  $x=\frac{n}{y}$ and consider  
the function $$G_n(y) := S_n \left(\frac{n}{y}\right) \in \mathbb{C}(y).$$ Then for every 
fixed $n\in \mathbb{N}$ all the zeros of $G_n$ belong to the set $$C_n :=\left \{z \in\mathbb{C}:\  
 \left|z + i\frac{n}{2b}\right| \geq  \frac{n}{2b},\    \left|z - i\frac{n}{2b}\right| \geq  \frac{n}{2b}\right \}.$$ 
 The sequence $G_n(y)$ converges 
uniformly on the compact sets to the entire function $$f(y) := \sum_{j=l}^m 
a_j e^{-j \lambda y} = Q (e^{- \lambda y})$$ as $n\to \infty.$ Each zero of the limiting 
entire function $f$ is the accumulation point of a sequence of zeros of $G_n$. 
Obviously if a sequence $\{ z_k\}_{k\in \mathbb{N}}$ has a limit $y_0$ and  
for all $ k\in \mathbb{N}$ we have $  z_k \in C_k, $  then $y_0$ belongs to the 
real axis. In the same way as in the previous case of hyperbolicity preservers we conclude that $\mbox{Re}\    
\lambda =0$  and that any non-zero root of $Q$ belongs to the circle  $\{ z : |z| = 
1\}.$ 

So, in both cases we have  $\lambda = i \beta,\   \beta \in \mathbb{R}\setminus\{0\},$
and $Q(t) := \sum_{j=l}^m a_j t^{j} = t^l \prod_{k=1}^{m-l}(t- e^{i\theta_k}), $
$\theta_k \in \mathbb{R}$ for $k=1, 2, \ldots , m-l.$ Then our linear operator $T$ has 
the following representation
\begin{equation}
\label{e7} T=  S_{i \beta}^l \prod_{k=1}^{m-l}(S_{i \beta}- e^{i\theta_k}I).
\end{equation}
Proof of the condition 2 in theorems \ref{th:mth1} and \ref{th:mth2} is based on the 
following fact on polynomials having all their roots on a horizontal straight line.\\  

{\bf Lemma}. {\it Let $T=S_{i\beta}-e^{i\theta}\cdot I ,$ where $\beta , \theta  \in \mathbb{R}.$ 
Suppose $P\in \mathbb{C}[z],$ and all the zeros of $P$ lie on a straight line $\{z: \ \ Im \ z=c\}.$ Then 
all the zeros of the polynomial $T(P)$ lie on the straight line $\{z: \ \ Im\ z=c + \beta/2 \}.$ 
} \\

{\it Proof of Lemma.}
Suppose $P\in\mathbb{C}[z]$ is an arbitrary polynomial having all zeros on the straight
line $\{ z: \ \mbox{Im}\  z =c  \}$, that is $$P(z)=C\prod_{j=1}^n (z-d_j -c i),$$ where $C\ne 0,
d_j, c\in\mathbb{R}.$ Let us investigate the possible zero location of $(S_{i \beta}- e^{i\theta}I)(P).$
We have $$(S_{i \beta}- e^{i\theta}I)(P)(z_0)=0 \ \Leftrightarrow \  
\prod_{j=1}^n\frac{z_0-d_j -c i-\beta i}{z_0-d_j -c i} = e^{i\theta}.$$ 
We observe that for all $j=1, 2, \ldots , n,$ the following is true
 $$\left |\frac{z-d_j -c i-\beta i}{z-d_j -c i}\right | <1 \ \mbox{whenever} \ \mbox{Im}\  z  > c+ \beta/2$$  
and $$\left |\frac{z-d_j -c i-\beta i}{z-d_j -c i}\right | >1 \ \mbox{whenever} \ \mbox{Im}\  z  <
 c+ \beta/2.$$  Thus, all the zeros of  $(S_{i \beta}- e^{i\theta}I)(P)$ belong to the  
line $\{ z: \ \mbox{Im}\  z =c +\beta/2 \}$ provided that all the zeros  of $P$ belong to the line 
$\{ z: \ \mbox{Im}\  z =c  \}.$ Lemma is proved. $\Box$ \\

Let's prove the necessity of the condition 2 in Theorem \ref{th:mth1}.
If $P$ is a hyperbolic polynomial, then all its zeros belong to the line $\{z: \ \ Im \ z=0\}.$ It follows
 from the above lemma that  all the zeros of  $\left(\prod_{k=1}^{m-l}
(S_{i \beta}- e^{i\theta_k}I)\right)(P) $  are on the line $\{ z: \ \mbox{Im}\  z =(m-l)\beta/2  \}.$ 
By virtue of (\ref{e7})  all the zeros of $ T(P) $  are on the line  $ \{ z: \ \mbox{Im}\  z =
(m-l)\beta/2 +l\beta \}.$ Therefore,  if all the zeros of 
$T(P)$ are real for every hyperbolic polynomial $P$,  then  $l = -m.$\\

Now we prove the necessity of the condition 2 in theorem \ref{th:mth2}, that is for strip preservers. 
Suppose that a linear operator of the form (\ref{e7}) preserves 
the set of complex polynomials having all zeros in the strip $\Pi_b :=\{ z:\  |\mbox{Im}\  z|\leq b \}.$
Consider any polynomial $P$ with all its zeros on the line $\{ z: \ \mbox{Im}\  z = b  \}.$ 
By the above lemma all zeros of $T(P)$ belong to the line $$\{ z: \ \mbox{Im}\  z = b +(m-l)\beta/2 +l\beta \}.$$
If a polynomial $P$ has all its zeros on the line $\{ z: \ \mbox{Im}\  z = - b  \},$ then
all the zeros of $T(P)$ are on the line $$\{ z: \ \mbox{Im}\  z = -b +(m-l)\beta/2 +l\beta \}.$$
Since the operator $T$ preserves the strip $\Pi_b$ this is possible only if
$$|\pm b  +(m-l)\beta/2 +l\beta| \leq b.$$
Thus $$(m-l)\beta/2 +l\beta = 0 \ \ \Leftrightarrow \ \l = -m.$$
That completes proof of the necessity of the conditions 1,2,3
in theorems \ref{th:mth1} and \ref{th:mth2}.\\
Note that the condition 4 in theorem \ref{th:mth1} provides the fact that for any polynomial $P \in \mathcal{HP}$  
the coefficients of the polynomial $T(P)$ are real.\\

The sufficiency of the conditions 1,2,3,4 in theorem \ref{th:mth1} follows from Theorem A and the 
remark that any hyperbolicity preserver $T$ is a composition of  linear operators of 
the form (\ref{e33}).\\

Let us prove the sufficiency of conditions 1, 2, 3 in theorem \ref{th:mth2}, that is for a strip preservers. Suppose that a polynomial 
$P$ has all zeros in the strip 
$\Pi_b :$ 
 $$P(z)=C\prod_{j=1}^n (z-z_j),$$ where $C\ne 0, |\mbox{Im}\  z_j  |\leq b,$ 
$j=1, 2, \ldots,n.$ 

Let us investigate the possible zero location of $(S_{i \beta}- e^{i\theta}I)(P).$
We have $$(S_{i \beta}- e^{i\theta}I)(P)(z_0)=0  \ \Leftrightarrow \ \prod_{j=1}^n 
\frac{z_0-z_j-\beta i}{z_0-z_j} = e^{i\theta}.$$ 
We observe that for all $j=1, 2, \ldots , n,$ the following is true:
 $$\left |\frac{z_0-z_j-\beta i}{z_0-z_j} \right | <1 \ \mbox{whenever} \ \mbox{Im}\  z  > b + \beta/2$$  
and $$\left |\frac{z_0-z_j-\beta i}{z_0-z_j}\right | >1 \ \mbox{whenever} \ \mbox{Im}\  z  < -b + \beta/2.$$

Hence all the zeros of  $(S_{i \beta}- e^{i\theta}I)(P)$ belong to the strip  
$\{ z:\  |\mbox{Im}\  z -\beta/2 |\leq b \},$ and all the zeros of a linear operator $T$ of the form (\ref{e7}) lie in the strip
$$\{ z:\  |\mbox{Im}\  z -(m+l)\beta/2 |\leq b \}.$$
Thus a linear operator $T$ of the form (\ref{e7}) under condition  $l = -m$  preserves the set of complex polynomials having all zeros in the strip 
$$\Pi_b \ = \ \{ z:\  |\mbox{Im}\  z|\leq b \}.$$  Theorems \ref{th:mth1} and \ref{th:mth2} are proved. $\Box$

\setcounter{equation}{0}
\section{Proof of Theorem  \ref{th:mth3}  }\label{section:proof.thm1.3}

{\it Proof of Theorem\ref{th:mth3}.}  For every $h>0, \theta \in \mathbb{R},$ and  every $f \in \mathcal{L-P}$  all 
the zeros of $T_{\theta, h}(f)$  are real since $f$ is the uniform limit, on compact subsets of
$ \mathbb{C}$, of polynomials with only real zeros, and, as it was mentioned before, $T_{\theta, h}:
\mathcal{HP} \to \mathcal{HP}.$ We need to prove only the simplicity of zeros of $T_{\theta, h}(f).$

Let $f \in \mathcal{L-P}$ have the representation (\ref{e2}), $f \not\equiv 0,$ and suppose that $x_0 \in \mathbb{R}$ is
the multiple root of $g(z) := T_{\theta, h}(f)(z).$ Then $g(x_0)=0,\  g^{\prime} (x_0)=0,$ or
$$  e^{i\theta} f(x_0 +ih) = e^{-i\theta} f(x_0 -ih),\quad   e^{i\theta} f^{\prime}(x_0 +ih) = e^{-i\theta} f^{\prime}(x_0 -ih), $$
whence $\frac{f^{\prime}}{f}(x_0 +i h) = \frac{f^{\prime}}{f}(x_0 -i h).$ By (\ref{e2}) we have
$$ \frac{f^{\prime}}{f}(z) = \frac{n}{z} -2az +b + \sum_{k=1}^\infty \frac{z}{x_k(z-x_k)}, $$
hence we obtain
$$  \frac{n}{x_0 +i h} -2a(x_0 +i h) +b + \sum_{k=1}^\infty \frac{x_0 +i h}{x_k(x_0 +i h-x_k)}=
  \frac{n}{x_0 -i h} -2a(x_0 -i h) +b + \sum_{k=1}^\infty \frac{x_0 -i h}{x_k(x_0 -i h-x_k)}, $$
or
$$  \frac{n}{x_0 +i h} -2aih + \sum_{k=1}^\infty \frac{x_0 +i h}{x_k(x_0 +i h-x_k)}=
 \frac{n}{x_0 -i h} +2aih+ \sum_{k=1}^\infty \frac{x_0 -i h}{x_k(x_0 -i h-x_k)}.$$
We compare the imaginary parts of the left hand and right hand sides:
$$\frac{-nh}{x_0^2 +h^2}- 2ah - \sum_{k=1}^\infty \frac{x_k h}{x_k((x_0 -x_k)^2 + h^2)}=
\frac{nh}{x_0^2 +h^2} + 2ah + \sum_{k=1}^\infty \frac{x_k h}{x_k((x_0 -x_k)^2 + h^2)}. $$
Since $h\ne 0$ we conclude that
$$\frac{n}{x_0^2 +h^2} + 2a + \sum_{k=1}^\infty \frac{1 }{(x_0 -x_k)^2 + h^2} = 0. $$
But $n\geq 0, a\geq 0, (x_0-x_k)^2 \geq 0,$ whence we get that $f$ is a constant function,
so $T_{\theta, h}(f)$ is a constant function and we are done. $\Box$

\setcounter{equation}{0}
\section{The Walsh convolution, Proof of Theorems \ref{th:mth5} and \ref{th:mth4} } \label{section:proof.thm1.4}

\begin{definition}[see, for example, {\cite[Chapter~5,~\textsection 3,~Problem~139]{PS}}]
\label{th:r1}
{\it Two complex polynomials $P$ and $Q$ of degree $n$ are called apolar if }
\begin{equation}
\sum_{k=0}^{n}(-1)^k \ P^{(k)}(0)\cdot Q^{(n-k)}(0)=0.
\end{equation}
\end{definition} 
 
The following famous theorem due to J.H. Grace states that the complex zeros of two apolar polynomials 
cannot be separated by a straight line or by a circle.\\ 

{\bf Theorem D} (J.H. Grace, see, for example,  \cite[Chapter~5,~\textsection 3,~Problem~145]{PS}).    
Suppose $P$ and $Q$ are two apolar polynomials of degree $n \geq 1.$ If all zeros of $P$ lie 
in a circular region $C,$ then $Q$ has at least one zero in $C.$ (A circular region is a closed or 
open half-plane, disk or exterior of a disk).\\

\begin{definition} [\cite{Walsh}]
\label{th:r2}
{\it For any two complex polynomials $P$ and $Q$ of degree $n$ the Walsh convolution is defined as follows }
\begin{equation}
P\boxplus Q \ (x) \ = \ \sum_{k=0}^{n} P^{(k)}(0)\cdot Q^{(n-k)}(x).
\end{equation}
\end{definition} 

By comparing formulas (\ref{th:r1}) and (\ref{th:r2}) we observe that
\begin{equation}
\label{th:r3}
P\boxplus Q \ (x_0)=0 \ \Leftrightarrow \ P(-x) \ \mbox{and} \ Q(x+x_0) \  \mbox{are apolar.}
\end{equation}

The following well-known fact was probably first proved by J.L.\,Walsh (\cite{Walsh}).

{\bf Theorem D}. {\it 1. For any two hyperbolic polynomials $P$ and $Q$ of degree $n,$ their Walsh convolution $P\boxplus Q$ is also hyperbolic.

2. If in addition all zeros of the polynomial $P$ lie in the interval $[\alpha, \ \beta]$, and all zeros of the polynomial $Q$ lie in the interval $[\gamma, \ \delta]$, then all zeros of the polynomial $P\boxplus Q$ lie in the interval $[\alpha+\gamma, \ \beta+\delta].$ }

For a reader's convenience we provide a proof of this theorem.

{\it Proof of Theorem D.} Let's prove the first statement. Assume that $P\boxplus Q(x_0)=0,$ but $ \mbox{Im}\ x_0=b \neq 0.$ 
So, by (\ref{th:r3}) the polynomial $P(-x)$ and $Q(x+x_0)$ are apolar.  Using the fact that the polynomials $P$ and $Q$ are 
hyperbolic, we conclude that all zeros of  $ P(-x)$  belong to the line $\{\mbox{Im}\ z=0\},$ while all zeros of $ Q(x+x_0)$  
belong to the line $ \{\mbox{Im}\ z=-b\}.$  Hence the zeros of the polynomials $P(-x)$ and $Q(x+x_0)$ can be separated 
by a straight line, which contradicts to the Grace's theorem.
The first statement of Theorem D is proved.\\

Now we prove the second statement. Consider any root $x_0$ of $P\boxplus Q.$ Given all zeros of $P(-x)$ lie in the 
interval $[-\beta, \ -\alpha]$ and all zeros of $Q(x+x_0)$ lie in the interval $[\gamma-x_0, \ \delta-x_0],$  the Grace's theorem 
provides an existence of a point $\zeta \in \mathbb{R}$ such that 
\begin{equation}
\label{th:r4}
-\beta \ \leq \ \zeta \ \leq \ -\alpha \ \ \mbox{and} \ \ \gamma-x_0 \ \leq \ \zeta \ \leq \ \delta-x_0.
\end{equation}
Thus, $$ \alpha+\gamma \ \leq \ x_0 \ \leq \ \beta+\delta \ .$$
Theorem D is proved. $\Box$\\

{\it Proof of Theorem \ref{th:mth5}.} For any $\theta \in \mathbb{R}$ and $h>0$ we consider the operator (\ref{e5}).
Denote by
\begin {equation} \label{R2} G_n(x,\theta, h)= T_{\theta, h}(x^n) =\frac{e^{i\theta}(x+ih)^n-e^{-i\theta}(x-ih)^n}{i}, \ \mbox{if}\  \sin\theta \neq 0,\ \end{equation}
and 
\begin {equation} \label{R3} G_n(x,2\pi m, h)=- G_n(x,\pi+2\pi m,h)=G_n(x,0,h)=T_{0, h}(x^n)=\frac{(x+ih)^n-(x-ih)^n}{i},  
\end{equation}
where $\ m\in\mathbb{Z}, \  n=0,1,2, \ldots.$

Comparing these formulas with (\ref{E2}) and (\ref{E3}) we observe that for
 every $\theta \in \mathbb{R},$  $h>0$ and $n\in \mathbb(N):$
\begin {equation} \label{R4}
G_n(x, \theta, h)=h^n Q_n\left(\frac{x}{h}, \theta \right). 
\end{equation}

So, for every $n \in \mathbb{N}$ the polynomial $G_n(x, \theta, h)$ is hyperbolic with the zeros 
$h\cdot x_k, \ k=1,2,\ldots, N,$ where $x_k$ are roots of the polynomial $Q_n(x, \theta),$ which together with the number $N$ are described by (\ref{a3}).\\

Let's apply the operator $T_{\theta, h}$
to a hyperbolic polynomial $$P(x)=\sum_{k=0}^n \frac{1}{k!}P^{(k)}(0) x^k.$$ We obtain

$$T_{\theta, h}(P)\ (x) \  = \ \sum_{k=0}^n \frac{1}{k!}\ P^{(k)}(0) \  T_{\theta, h}(x^k) \ =\ 
\sum_{k=0}^n \frac{1}{k!}\ P^{(k)}(0)\cdot G_k(x,  \theta, h).$$

It is easy to show that $$G_k(x, \theta, h)\ =\ \frac{G_n^{(n-k)}(x, \theta, h)} {n(n-1)(n-2)\cdots(k+1)}\ =\ 
\frac{k!}{n!}\ G_n^{(n-k)}(x, \theta, h). $$
Thus
\begin{equation}
\label{th:r4}
T_{\theta, h}(P)\ (x) =\ \frac{1}{n!}\ \sum_{k=0}^n \ P^{(k)}(0) \  G_n^{(n-k)}(x, \theta,h)\ =\ 
\frac{1}{n!}\ P \boxplus G_n(x, \theta, h ).
\end{equation}

Now the statement of Theorem \ref{th:mth5} follows immediately from the second statement of Theorem D and (\ref{R4}). $\Box$\\

 Let us fix a polynomial $P\in \mathbb{C}[x].$
One can consider a linear operator $T_{P}: \ \mathbb{C}[x] \ \to \ \mathbb{C}[x]$ acting as follows: $$  \ T_{P}(Q) \ = \ P \boxplus Q.$$
If $P\in \mathcal{HP},$ then according to Theorem D the operator $T_{P}$ is a hyperbolicity preserver. Obviously, $T_{P}$ commutes with any shift operator. It follows from Theorem C that $T_{P}$ does not decrease mesh. Using the commutative property of the Walsh convolution (see, for example \cite{Walsh}) we obtain the following result proved in \cite{BKS}.\\

{\bf Theorem E} (\cite{BKS}). {\it Let $P$ and $Q$ be hyperbolic polynomials of degree $n$. Then \begin{equation}
\label{th:r5} \mbox{mesh}\ (P\boxplus Q) \ \geq \ \max \ (\mbox{mesh}\ (P),\ \mbox{mesh}\ (Q)).
\end{equation}
}

{\it Proof of Theorem \ref{th:mth4}.} Let us apply the operator $T_{\theta, h}$ to a hyperbolic polynomial $P.$ It follows from 
(\ref{th:r4}) and (\ref{th:r5}) that 
$$\mbox{mesh}\ T_{\theta, h}(P)(x)\geq 
\mbox{mesh}\ (P\boxplus G_n(x, \theta, h)) \ \geq \ \max \ (\mbox{mesh}\ (P(x)),\ \mbox{mesh}\ ( G_n(x, \theta, h))).$$

Since $G_n(x, \theta, h)= T_{\theta, h}(x^n)$ this completes the proof of Theorem \ref{th:mth4} as well as presents another vision of the fact that the operator $T_{\theta, h}$ doesn't decrease the mesh of a hyperbolic polynomial. $\Box$\\

As we mentioned in Introduction of this paper Theorem C was established by S.Fisk (\cite[p. 226, Lemma 8.25]{Fisk}). However, it was made not in a very lucid way. We provide a proof of Theorem C. Our proof is based on the following fact known under the name Obreschkov's theorem (see, for example \cite[p. 10]{O}, although it has been
rediscovered many times by different authors in the past). See  \cite{KSV} for the analogous proof for the minimal quotent of roots 
instead of the minimal distance.\\

{\bf Theorem F} (N.\,Obreschkov, \cite[p. 10]{O}, or \cite[p. 10, Proposition 1.35]{Fisk}). Given two real polynomials $P$ and $Q$ of the same degree one
has that the pencil $c P(x)+d Q(x), \ c, \ d \ \in \mathbb{R},$ consists of hyperbolic polynomials if and
only if $P$ and $Q$ have  all real and (non-strictly) interlacing roots.\\
 
 {\it Proof of Theorem C.} First, note that for any hyperbolic polynomial $P$ the zeros of $P(x)$ and $P(x+\lambda)$ are 
 (non-strictly) interlacing if and only if $\lambda \leq mesh(P)$.

Assume that a linear operator $A: \ \mathcal{HP} \ \to \ \mathcal{HP}$ commutes with any shift operator $S_b, \ b \in 
\mathbb{R},$ but there exists a hyperbolic polynomial $P$ such that $\mbox{mesh}\ (A(P)) \ < \ \mbox{mesh}\ (P).$ It means that we can find such a real number $\lambda$ that 
$$\mbox{mesh}\ (A(P)) \ < \ \lambda \ < \ \mbox{mesh}\ (P).$$  So, as we mentioned above the roots of the hyperbolic polynomials $P(x)$ and $P(x+\lambda)$ are interlacing, while the roots of the hyperbolic polynomials $A(P)(x)$ and $A(P)(x+\lambda)$ are not. It follows from 
Theorem F that there are such two numbers $c,\ d \in \mathbb{R}$ that \begin{equation}
\label{th:r6}
 c A( P)(x)+d A(P)(x+\lambda) \notin \mathcal{HP},
\end{equation}
while
\begin{equation}
\label{th:r7}
c P(x)+d P(x+\lambda) \in \mathcal{HP}.  
\end{equation}

Since the operator $A$ is a hyperbolicity preserver, by virtue of (\ref{th:r6}) we have
$$ A\left( c P(x)+d P(x+\lambda)\right)\ = \ c A(P)(x)+d A(S_{\lambda}(P))(x) \in \mathcal{HP}. $$
On the other hand since the operator $A$ commutes with a shift operator,  the following is true
$$\mathcal{HP} \ni \ c A(P)(x)+d A(S_{\lambda}(P))(x)\  = \  c A(P)(x)+d S_{\lambda}(A(P))(x) \ = \ c A( P)(x)+d A(P)(x+\lambda),$$
which contradicts to (\ref{th:r7}).\\
Theorem C is proved. $\Box$

\section{Proof of Theorem \ref{th:mth6} }\label{section:proof.thm1.11}

Since $D_n(x, 2\pi k,h)=-D_n(x,\pi+ 2\pi k,h)=D_n(x, 0,h)$ we will consider only the cases: $sin\theta \ne 0$ and $ \theta=0.$

The following fact about properties of the polynomials $Q_n(x,\theta)$ is obvious. \\

{\bf Statement.} For each $n=2,3,\ldots$ the following relations are true
 \begin {equation} Q_n^{'}(x,\theta)=nQ_{n-1}(x,\theta),  \ \theta \in \mathbb{R}; \label{a4} \end{equation}
\begin {equation} Q_n^{'}(x_j,\theta)=2 \sin\theta \prod_{k\neq j} (x_j-x_k ), \mbox{if} \  sin\theta \ne 0, \ \mbox{and} \ \ Q_n^{'}(x_j,0)=2 n\prod_{k\neq j} (x_j-x_k ); \label{a5} \end{equation}
\begin {equation} Q_{n-1}(x_j,\theta)=\frac{2 \sin\theta}{n}\prod_{k\neq j} (x_j-x_k ), \mbox{if} \  sin\theta \ne 0, \ \mbox{and} \ \ Q_{n-1}(x_j,0)=2 \prod_{k\neq j} (x_j-x_k ). \label{a6} \end{equation}

Let's divide $D_n(x,\theta,h)$ by $h^n$  $$\frac{1}{h^n}D_n(x, \theta,h )=\frac{1}{i}\left\{e^{i\theta}\left(\frac{x}{h}+i\right)^n-e^{-i\theta}\left(\frac{x}{h}-i\right)^n\right\}+\frac{a}{ih}\left\{e^{i\theta}\left(\frac{x}{h}+i\right)^{n-1}-e^{-i\theta}\left(\frac{x}{h}-i\right)^{n-1}\right\}$$ $$+\frac{b}{ih^2}\left\{e^{i\theta}\left(\frac{x}{h}+i\right)^{n-2}-e^{-i\theta}\left(\frac{x}{h}-i\right)^{n-2}\right\}+\sum_{k=0}^{n-3}\frac{c_k}{ih^{n-k}}\left\{e^{i\theta}\left(\frac{x}{h}+i\right)^{k}-e^{-i\theta}\left(\frac{x}{h}-i\right)^{k}\right\}.$$

Denote by \begin{equation} t=\frac{x}{h}.\label {b1} \end {equation} Using the formulas from Example \ref{E1} we can reformulate our problem as follows: to describe asymptotic behavior of the roots $t_1(h,\theta),t_2(h,\theta),\ldots,t_{N}(h,\theta)$ (as before $N=n$ if $sin \theta \ne 0,$ and $N=n-1$ if $\theta = 0$) of the equation
\begin {equation} Q_n(t,\theta)+\frac{a}{h}Q_{n-1}(t,\theta)+\frac{b}{h^2}Q_{n-2}(t,\theta)+\sum_{k=0}^{n-3}\frac{c_k}{h^{n-k}}Q_{k}(t,\theta)=0 \label {a7} \end {equation}
as $h \to \infty$  for each $\theta \in \mathbb{R}$.\\
We fix an integer number $ j=1,2,\ldots, N.$ Denote by 
\begin {equation} P_n(t,\theta)=2 \sin\theta \prod_{k\neq j} (t-x_k ), \mbox{if} \  sin\theta \ne 0, \ \mbox{and} \ \ P_n(t,0)=2 n\prod_{k\neq j} (t-x_k ).\label{o1} \end{equation}
Since $Q_n(t,\theta)=(t-x_j) P_n(t,\theta),$ one can derive the following properties of $P_n$ from (\ref{a4}), (\ref{a5}) and (\ref{a6}):
\begin {equation} P_n(x_j, \theta)=Q_n^{'}(x_j,\theta)=n Q_{n-1}(x_j,\theta),  \label{o2} \end{equation}
\begin {equation} P_n^{'}(x_j, \theta)=\frac{1}{2}Q_n^{''}(x_j,\theta)=\frac{n(n-1)}{2}Q_{n-2}(x_j,\theta).   \label{o3} \end{equation}
By Hurwitz theorem there exists such a number $\rho>0$ that for big enough values of the number $h$ the circle $|t-x_j|<\rho$ contains only one root of the equation (\ref{a7}), and this root is $t_j(h,\theta)$, that is 
\begin{equation} |t_j(h,\theta)-x_j|<\rho,\quad |t_k(h,\theta)-x_j|\geq \rho, \ \ k\neq j. 
\label {c1} \end {equation}  
Therefore $P_n(t_j(h,\theta), \theta)\ne 0,$ and we can divide (\ref{a7}) by it. Thus the root $t_j(h,\theta)$ satisfies the equation 
$$ (t_j(h,\theta)-x_j)+\frac{a}{h} \frac{Q_{n-1}(t_j(h,\theta),\theta)}{P_n(t_j(h,\theta),\theta)}+\frac{b}{h^2}\frac{Q_{n-2}(t_j(h,\theta),\theta)}{P_n(t_j(h,\theta),\theta)}$$  \begin {equation}+\sum_{k=0}^{n-3}\frac{c_k}{h^{n-k}}\frac{Q_{k}(t_j(h,\theta),\theta)}{P_n(t_j(h,\theta),\theta)} =0. \label{a9} \end{equation} 
It follows from this equation and (\ref{c1}) that \begin{equation} t_j(h,\theta)=x_j +O\left(\frac{1}{h}\right). \label {a13} \end{equation}
Using the Taylor expansion formulas for the functions $\frac{Q_{n-1}(t,\theta)}{P_n(t,\theta)}$ and $\frac{Q_{n-2}(t,\theta)}{P_n(t,\theta)}$ about $x_j,$ and (\ref{o2}), (\ref{o3}), (\ref{a4}) we obtain

$$ \frac{Q_{n-1}(t,\theta)}{P_n(t,\theta)}=\frac{Q_{n-1}(x_j,\theta)}{P_n(x_j,\theta)}+$$ 
$$\frac{Q^{'}_{n-1}(x_j,\theta)P_n(x_j,\theta)-Q_{n-1}(x_j,\theta)P^{'}_n(x_j,\theta)}{P^2_n(x_j,\theta)}(t-x_j) +O\left((t-x_j)^2\right)=$$ \begin {equation}\frac{1}{n}+\frac{(n-1)}{2n}\cdot\frac{Q_{n-2}(x_j,\theta)}{Q_{n-1}(x_j,\theta)} (t-x_j) +O\left((t-x_j)^2\right), \label{a10}\end{equation} 
and \begin {equation} \frac{Q_{n-2}(t,\theta)}{P_n(t,\theta)}=\frac{Q_{n-2}(x_j,\theta)}{nQ_{n-1}(x_j,\theta)}+O\left((t-x_j)\right). \label{a11}\end{equation} 
The relations (\ref{a10}) and (\ref{a11}) allow us to rewrite (\ref{a9}) in the following way
$$ (t_j(h,\theta)-x_j)+\frac{1}{h}\cdot\frac{a}{n}+\frac{1}{h}\cdot\frac{a(n-1)}{2n}\cdot\frac{Q_{n-2}(x_j,\theta)}{Q_{n-1}(x_j,\theta)} \left(t_j(h,\theta)-x_j\right) +\frac{1}{h^2}\cdot \frac{b}{n}\frac{Q_{n-2}(x_j,\theta)}{Q_{n-1}(x_j,\theta)}$$ $$+\frac{a}{h} \cdot O\left((t_j(h,\theta)-x_j)^2\right)+\frac{b}{h^2}\cdot O\left((t_j(h,\theta)-x_j)\right)+$$ \begin {equation}\sum_{k=0}^{n-3}\frac{c_k}{h^{n-k}}\frac{Q_{k}(t_j(h,\theta),\theta)(t_j(h,\theta)-x_j)}{Q_n(t_j(h,\theta),\theta)} =0. \label{a12} \end{equation} 
 Let's put 
\begin{equation} t_j(h,\theta )-x_j=-\frac{a}{nh}+\left(\frac{a^2(n-1)}{2n^2h^2} -\frac{b}{nh^2}\right)\frac{Q_{n-2}(x_j,\theta )}{Q_{n-1}(x_j,\theta)} +\omega(h,\theta)\label {b2}, \end {equation}
and estimate the function $\omega (h,\theta)$ for big enough values of $h.$\\
Using (\ref{a13}) and (\ref{c1}), from (\ref{a12}) we obtain the following estimation 
\begin{equation} |\omega (h,\theta) |\leq \frac{K}{h^3},\label {b3} \end {equation} where $K$ is a constant.
Thus by virtue of (\ref{b1}) we have \\
$$X_j(h,\theta)=h\cdot t_j(h,\theta)= x_j\cdot h-\frac{a}{n}+\left( \frac{a^2(n-1)}{2n^2}  - 
\frac{b}{n}\right)\frac{Q_{n-2}(x_j,\theta)}{Q_{n-1}(x_j,\theta)}\cdot \frac{1}{h} +
O\left(\frac{1}{h^2}\right), \quad h\to \infty .$$
Theorem \ref{th:mth4} is proved.

In the same way we can obtain the more precise asymptotic formula:
$$X_j(h,\theta)=x_j\cdot h-\frac{a}{n}+\left( \frac{a^2(n-1)}{2n^2}  -\frac{b}{n}\right)\frac{Q_{n-2}(x_j,\theta)}{Q_{n-1}(x_j,\theta)}\cdot \frac{1}{h} $$
$$+\left(-\frac{a^3(n-1)(n-2)}{3n^3}+\frac{ab(n-2)}{n^2}-\frac{c}{n}\right)\frac{Q_{n-3}(x_j,\theta)}{Q_{n-1}(x_j,\theta)}\cdot \frac{1}{h^2} + O\left(\frac{1}{h^3}\right), \quad h\to \infty . $$
For this purpose we have to replace formulas (\ref{a10}) and (\ref{a11}) by the more accurate formulas
$$ \frac{Q_{n-1}(t,\theta)(t-x_j)}{Q_n(t,\theta)}=\frac{1}{n}+\frac{(n-1)Q_{n-2}(x_j,\theta)}{2nQ_{n-1}(x_j,\theta)} (t-x_j)+ $$  \begin {equation}\frac{1}{2}\left( \frac{2(n-1)(n-2)Q_{n-3}(x_j,\theta)}{3nQ_{n-1}(x_j,\theta)}-\frac{(n-1)^2 Q^2_{n-2}(x_j,\theta)}{2nQ^2_{n-1}(x_j,\theta)}\right)\cdot (t-x_j)^2+O\left((t-x_j)^3\right), \label{q1}\end{equation} 
and $$ \frac{Q_{n-2}(t,\theta)(t-x_j)}{Q_n(t,\theta)}=\frac{Q_{n-2}(x_j,\theta)}{nQ_{n-1}(x_j,\theta)}+$$
\begin {equation} \left( \frac{(n-2)Q_{n-3}(x_j,\theta)}{nQ_{n-1}(x_j,\theta)}-\frac{(n-1) Q^2_{n-2}(x_j,\theta)}{2nQ^2_{n-1}(x_j,\theta)}\right)\cdot (t-x_j)+O\left((t-x_j)^2\right). \label{q2}\end{equation}
Additionally, we write down the first term of the Taylor expansion for the function $\frac{Q_{n-3}(t,\theta)}{Q_n(t,\theta)}(t-x_j):$ 
\begin {equation} \frac{Q_{n-3}(t,\theta)(t-x_j)}{Q_n(t,\theta)}=\frac{Q_{n-3}(x_j,\theta)}{nQ_{n-1}(x_j,\theta)}+O\left((t-x_j)\right). \label{q3}\end{equation} 
After that concidering the formulas  (\ref{q1}), (\ref{q2}) and (\ref{q3}) we make corresponding changes in the formula (\ref{a12}).  Putting 
$$ t_j(h,\theta )-x_j=-\frac{a}{nh}+\left(\frac{a^2(n-1)}{2n^2h^2} -\frac{b}{nh^2}\right)\frac{Q_{n-2}(x_j,\theta )}{Q_{n-1}(x_j,\theta)} +$$ \begin{equation} \left(-\frac{a^3(n-1)(n-2)}{3n^3h^3}+\frac{ab(n-2)}{n^2h^3}-\frac{c}{nh^3}\right)\frac{Q_{n-3}(x_j,\theta)}{Q_{n-1}(x_j,\theta)} +\omega(h,\theta)\label {q4}, \end {equation}
we obtain the following estimation for the function $\omega (h,\theta)$
\begin{equation} |\omega (h,\theta) |\leq \frac{K}{h^4},\label {q5} \end {equation} where $K$ is a constant. $\Box$

\end{document}